\documentclass[a4paper,notitlepage]{article}

\usepackage[latin1]{inputenc} 
\usepackage[T1]{fontenc}
\usepackage{amsmath}
\usepackage{amsthm}
\usepackage{amssymb}   
\usepackage{mathrsfs}  

\theoremstyle{plain}
\newtheorem{thm}{Theorem}[section]
\newtheorem{lemma}[thm]{Lemma}
\newtheorem{prop}[thm]{Proposition}

\theoremstyle{definition}
\newtheorem{define}{Definition}[section]
\newtheorem*{define*}{Definition}

\theoremstyle{remark}
\newtheorem{remark}[thm]{Remark}

\newcommand{\zmk}{\ensuremath{\zeta_m^{(k)}}}

\newcommand{\reals}{\ensuremath{\mathbb{R}}}

\newcommand{\C}{\ensuremath{\mathbb{C}}}
\newcommand{\BB}{\ensuremath{\left[\begin {array}{ccc} e^{it} & 0\\
\noalign{\medskip} 0 &e^{-it}\end {array}\right] }}
\newcommand{\B}{\ensuremath{d \omega \left(\begin {array}{ccc} i & 0\\
\noalign{\medskip} 0 &-i\end {array}\right) f}}
\title{Local Riemann Hypothesis for complex numbers}
\author{Rikard Olofsson}

\begin{document}
\maketitle
\begin{abstract}
In this paper a special class of local $\zeta$-functions is studied. The main 
theorem states that the functions have all zeros on the line $\Re (s)=1/2.$
This is a natural generalization of the result of Bump and Ng stating that the 
zeros of the Mellin transform of Hermite functions have $\Re (s)=1/2.$
\end{abstract}

\section{Introduction}
In the study of Hecke L-functions, Tate \cite{Tate},\cite{Tate2} 
defined local $\zeta$-functions 
\begin{align*}
\zeta(s,\nu ,f)=\int_{F^{\times}}f(x)\nu (x)|x|^sd^{\times}x,
\end{align*}
where $F$ is a local field, $f$
is a Schwartz function of $F,$ $\nu$ is a character of $F^{\times}$ and 
integration is taken with respect to Haar measure on $F^{\times}.$ 
Weil \cite{Weil} 
introduced a representation $\omega=\omega_{\psi}$ of the metaplectic group 
$\widetilde{SL}(2,F)$ for each nontrivial additive character 
$\psi$ of $F$.
The {\it Local Riemann Hypothesis (LRH)}, as formulated in \cite{lrh},
is the assertion that if $f$ is taken from some 
irreducible invariant subspace of the restriction of this representation to 
a certain compact subgroup $H$ of $\widetilde{SL}(2,F),$ 
then in fact all zeros of  
$\zeta(s,\nu ,f)$ lie on the line $\Re(s)=1/2.$ The phenomenon was first 
observed by Bump and Ng and they proved that the zeros of the Mellin transform 
of Hermite functions lie on the line, this corresponds to LRH for 
$F=\mathbb{R}$ \cite{Ng}. 
LRH has also been proved for
$F$ having odd characteristics by Kurlberg \cite{Kurl} 
and disproved for $F=\C$ by 
Kurlberg \cite{Kurl}. In all cases above $H$ is the unique maximal compact 
subgroup of $SO(2,F),$
for $F=\mathbb{R}$ and for $F$ with characteristic congruent to $3$ modulo $4,$
$H$ is nothing but $SO(2,F),$ since this already is compact.
In \cite{lrh} Bump, Choi, Kurlberg and Vaaler offer generalizations of LRH to 
higher dimensions along with 
two different proofs of the case $F=\mathbb{R}$ and $H=SO(2).$ In this paper 
we prove: 

\begin{thm}
\label{lrh}
If $f$ belongs to an irreducible invariant subspace of the Weil 
representation restricted to $SU(2,\C)$ and $\zeta(s,\nu,f)\not\equiv 0,$ then 
all zeros of $\zeta(s,\nu,f)$ lie on the line $\Re (s)=1/2.$
\end{thm}
In other words, we prove that a slightly modified version of  
LRH (namely taking $H=SU(2,\C)$ rather than a compact subgroup of $SO(2,\C)$) 
holds for $F=\C .$ 

\begin{remark}
From now on we will restrict ourselves to the case where the local field is 
$\C.$
\end{remark}

\section{Acknowledgments}
I would like to thank Pär Kurlberg for suggesting this problem to me and for
all his help and encouragement. 



\section{The Weil representation}
The Weil (or the metaplectic) representation is an action on 
$S(\C)=\{f(z);f(x+iy)=g(x,y)\in \mathscr{S}(\mathbb{R}^2)\},$
where $\mathscr{S}(\mathbb{R}^2)$ is the Schwartz space. 
We will often think of the elements of $S(\C)$ 
not as functions of the complex variable $z,$ but rather as functions 
of the two real variables $x,y$ satisfying $z=x+iy.$ In agreement with 
that we write $dz$ and this is nothing but $dxdy,$ the Lebesgue measure of 
$\mathbb{R}^2.$ Sometimes we will also use 
the notation $\langle f,g \rangle =\int_{\C}f(z)g(z)dz.$ 
Let the additive character on $\C$ be $\psi(z)=e^{i\pi\Re(z)}$
and introduce the Fourier transform
\begin{align*}
\hat f(z)=\int_{\C}f(z')\psi(2zz')dz'.
\end{align*}
With this normalization, we find that $\hat{\hat{f}}(z)=f(-z).$ 

\begin{remark}
As noted in \cite{lrh} there is no loss of generality in assuming that 
the additive character is $\psi(z)=e^{i\pi\Re(z)}$ if the objective only is to 
prove LRH. Changing character does not preserve the irreducible subspaces, 
but the zeros of the ``corresponding $\zeta$-functions'' are preserved.
\end{remark}

$\widetilde{SL}(2,\C),$ the metaplectic double cover of $SL(2,\C),$ splits and 
we have \\
$\widetilde{SL}(2,\C)\cong SL(2,\C)\times C_2.$
Using this identification we write 
\begin{align*}
\left[\begin {array}{ccc} a & b\\
\noalign{\medskip}c & d\end {array}\right]=
\left(\left(\begin {array}{ccc} a & b\\
\noalign{\medskip}c & d\end {array}\right) ,1\right).
\end{align*}
The restriction of 
the metaplectic representation to $SU(2,\C)$ can now be written as
\begin{align*}
\left(\omega\left[\begin {array}{ccc} \alpha & - \bar\beta\\
\noalign{\medskip}\beta & \bar\alpha\end {array}\right]f\right)(z)=
\frac{1}{|\beta |} \int_{\mathbb{C}} \psi \left( \frac{1}{\beta }
\left( \alpha z^2-2z z'+\bar\alpha {z'} ^2
\right) \right) f(z')dz'.
\end{align*}
However, it is much more convenient
to see how $\omega$ acts on the generators of 
$SL(2,\C).$ This is given by
\begin{align*}
\left( \omega \left[ \begin {array}{ccc} 1 & t\\
\noalign{\medskip} 0 & 1\end {array} \right] f \right)(z)=\psi(tz^2)f(z),
\end{align*}
\begin{align*}
\left( \omega \left[ \begin {array}{ccc} 0 & 1\\
\noalign{\medskip} -1 & 0\end {array} \right] f \right)(z)=\hat{f} (z),
\end{align*}
and
\begin{align*}
\left( \omega \left[ \begin {array}{ccc} \alpha & 0\\
\noalign{\medskip} 0 & \alpha^{-1}\end {array} \right] f \right)(z)
=|\alpha|f(\alpha z).
\end{align*}

\begin{remark}
When we write $|\alpha|$ we mean the ordinary absolute value of $\alpha,$
not the ``absolute value'' of an element in a local field used by Tate. 
\end{remark}


In order to find the invariant subspaces of the action of $SU(2,\C)$ 
we could of course 
just as well study the restriction to $\mathfrak{su}(2,\C)$ of the corresponding 
Lie algebra representation $d\omega:\mathfrak{sl}(2,\C)\to End(S(\C))$ defined by
\begin{align*}
((d\omega~X)f)(z)=\frac{d}{dt}\left(\omega\left(\widetilde{exp}(tX)\right)
f\right)(z)|_{t=0},
\end{align*}
where $\widetilde{exp}$ is the exponential map 
$\mathfrak{sl}(2,\C)\to SL(2,\C)$ lifted to a map  
$\widetilde{exp}:\mathfrak{sl}(2,\C)\to \widetilde{SL}(2,\C).$ 
Since a natural basis for $\mathfrak{su}(2,\C)$ is  
\begin{align*}
\left\{\left( \begin {array}{ccc} 0 & 1\\
\noalign{\medskip} -1 & 0\end {array} \right),
\left( \begin {array}{ccc} 0 & i\\
\noalign{\medskip} i & 0\end {array} \right),
\left( \begin {array}{ccc} i & 0\\
\noalign{\medskip} 0 & -i\end {array} \right)\right\},
\end{align*}
our first objective is to calculate how $d\omega$ acts on $S(\C)$ 
for these vectors. From the definitions we immediately get
\begin{align*}
\left( d \omega \left(\begin {array}{ccc} 0 & 1\\
\noalign{\medskip} 0 & 0\end {array}\right) f  \right)  & =
\frac{d}{dt} \left( \omega \left[ \begin {array}{ccc} 1 & t\\
\noalign{\medskip} 0 & 1\end {array} \right] f \right) \Bigg\vert_{t=0} = 
\frac{d}{dt}\psi\left( t (x+iy)^2 \right) f\Big\vert_{t=0}\\ & = 
\frac{d}{dt}e^{i\pi t \left(x^2-y^2\right)}f\Big\vert_{t=0}=i\pi 
\left(x^2-y^2\right)f
\end{align*}
and
\begin{align*}
\left( d \omega \left(\begin {array}{ccc} 0 & i\\
\noalign{\medskip} 0 & 0\end {array}\right) f  \right)  & =
\frac{d}{dt} \left( \omega \left[ \begin {array}{ccc} 1 & it\\
\noalign{\medskip} 0 & 1\end {array} \right] f \right) \Bigg\vert_{t=0} = 
\frac{d}{dt}\psi\left( it (x+iy)^2 \right) f\Big\vert_{t=0}\\ & = 
\frac{d}{dt}e^{-i2\pi txy }f\Big\vert_{t=0}=-i2\pi xy f.
\end{align*}
Introducing the notation $\mathfrak{F}$ for the operator taking $f$ to its 
Fourier transform $\hat f$ we see that 
\begin{align*}
d \omega \left( \begin {array}{ccc} 0 & 0\\
\noalign{\medskip} -1 & 0\end {array} \right) & =\left( \omega
\left[\begin {array}{ccc} 0 & 1\\
\noalign{\medskip} -1 & 0\end {array}\right]\right)^{-1}
\left(d \omega \left(\begin {array}{ccc} 0 & 1\\
\noalign{\medskip} 0 & 0\end {array}\right)\right)
\left( \omega \left[\begin {array}{ccc} 0 & 1\\
\noalign{\medskip} -1 & 0\end {array}\right]\right)\\ & =
\mathfrak{F}^{-1}i\pi\left(x^2-y^2\right)\mathfrak{F}=-\frac{i}{4\pi}
\left(\frac{\partial^2}{\partial x^2}-\frac{\partial^2}{\partial y^2}\right)
\end{align*}
and
\begin{align*}
d \omega \left( \begin {array}{ccc} 0 & 0\\
\noalign{\medskip} -i & 0\end {array} \right) & =\left( \omega
\left[\begin {array}{ccc} 0 & 1\\
\noalign{\medskip} -1 & 0\end {array}\right]\right)^{-1}
\left(d \omega \left(\begin {array}{ccc} 0 & i\\
\noalign{\medskip} 0 & 0\end {array}\right)\right)
\left( \omega \left[\begin {array}{ccc} 0 & 1\\
\noalign{\medskip} -1 & 0\end {array} \right] \right) \\ & =
\mathfrak{F}^{-1}\left(-i 2 \pi x y\right)\mathfrak{F} = -\frac{i}{2\pi}
\frac{\partial^2}{\partial x \partial y}.
\end{align*}
Hence we have that 
\begin{align*}
d \omega \left( \begin {array}{ccc} 0 & 1\\
\noalign{\medskip} -1 & 0\end {array} \right)=
d \omega \left( \begin {array}{ccc} 0 & 1\\
\noalign{\medskip} 0 & 0\end {array} \right)+
d \omega \left( \begin {array}{ccc} 0 & 0\\
\noalign{\medskip} -1 & 0\end {array} \right)=i\pi 
\left(x^2-y^2\right)-\frac{i}{4\pi}
\left(\frac{\partial^2}{\partial x^2}-\frac{\partial^2}{\partial y^2}\right)
\end{align*}
and
\begin{align*}
d \omega \left( \begin {array}{ccc} 0 & i\\
\noalign{\medskip} i & 0\end {array} \right)=
d \omega \left( \begin {array}{ccc} 0 & i\\
\noalign{\medskip} 0 & 0\end {array} \right)-
d \omega \left( \begin {array}{ccc} 0 & 0\\
\noalign{\medskip} -i & 0\end {array} \right)=-i2\pi xy+\frac{i}{2\pi}
\frac{\partial^2}{\partial x \partial y}.
\end{align*}
Finally we get that
\begin{align*}
\left(\B\right) & =
\frac{d}{dt} \left( \omega \BB f \right) (x+iy) \Bigg\vert_{t=0} = 
\frac{d}{dt} f(e^{it}(x+iy))\Big\vert_{t=0}\\ & = 
\frac{d}{dt} f(x\cos t - y\sin t + i(y\cos t +x \sin t))\Big\vert_{t=0}\\ &= 
-y\frac{\partial f}{\partial x} +x \frac{\partial f}{\partial y}.
\end{align*}

\begin{define}
Let $f_{m,n} (x+iy) = H_m(\sqrt{2\pi}x)H_n(\sqrt{2\pi}y)
e^{-\pi\left(x^2+y^2\right)}$ where 
$H_n(x)=(-1)^ne^{x^2}\frac{d^n}{dx^n}e^{-x^2}$ 
are the Hermite polynomials. 
\end{define}

\begin{prop}
$W_m=\bigoplus_{j=0}^m\mathbb{C}f_{j,m-j}$ are invariant subspaces of the Weil 
representation restricted to $\mathfrak{su}(2,\mathbb{C})$.
\end{prop}

\begin{proof}
We can write (see for instance \cite{Sze}) $f_{m,n} (x+iy)=h_m(x)h_n(y)$ where
$h_m$ satisfy
\begin{align*}
\left(x^2-\frac{1}{4\pi^2}\frac{d^2}{dx^2}\right)h_m=\frac{2m+1}{2\pi}h_m.
\end{align*}
Hence we have 
\begin{align*}
d \omega \left( \begin {array}{ccc} 0 & 1\\
\noalign{\medskip} -1 & 0\end {array} \right)f_{m,n} & =
\left(i\pi 
\left(x^2-y^2\right)-\frac{i}{4\pi}
\left(\frac{\partial^2}{\partial x^2}-\frac{\partial^2}{\partial y^2}\right)
\right)f_{m,n}\\ & =i\pi\left(\frac{2m+1}{2\pi}-\frac{2n+1}{2\pi}\right)f_{m,n}=
i(m-n)f_{m,n}.
\end{align*}
Using the recurrence formulas $H_{n+1}(x)=2xH_n(x)-2nH_{n-1}(x)$ and 
$H_n'(x)=2nH_{n-1}(x)$ \cite{Sze} we get
\begin{align*}
d \omega \left( \begin {array}{ccc} i & 0\\
\noalign{\medskip} 0 & -i\end {array} \right)f_{m,n} & =-y\frac{\partial f_{m,n}}
{\partial x} +x \frac{\partial f_{m,n}}{\partial y}\\ & =
-y\left(\sqrt{2\pi}2mf_{m-1,n}-2\pi xf_{m,n}\right)\\ &
+x\left(\sqrt{2\pi}2nf_{m,n-1}-2\pi yf_{m,n}\right)\\ &=
\sqrt{2\pi}\left(-2myf_{m-1,n}+2nxf_{m,n-1}\right)\\&=
-2m\frac{f_{m-1,n+1}+2nf_{m-1,n-1}}{2}
+2n\frac{f_{m+1,n-1}+2mf_{m-1,n-1}}{2}\\
&=nf_{m+1,n-1}-mf_{m-1,n+1}
\end{align*}and
\begin{align*}
d \omega \left( \begin {array}{ccc} 0 & i\\
\noalign{\medskip} i & 0\end {array} \right)f_{m,n}& 
=\frac{1}{2}d\omega 
\left[ \left( \begin {array}{ccc} i & 0\\
\noalign{\medskip} 0 & -i\end {array} \right), 
\left( \begin {array}{ccc} 0 & 1\\
\noalign
{\medskip} -1 & 0\end {array} \right)\right]f_{m,n}\\ & =
\frac{1}{2} d \omega \left( \begin {array}{ccc} i & 0\\
\noalign{\medskip} 0 & -i\end {array} \right)i(m-n)f_{m,n}\\
&-\frac{1}{2} d \omega \left( \begin {array}{ccc} 0 & 1\\
\noalign{\medskip} -1 & 0\end {array} \right)\left(nf_{m+1,n-1}-mf_{m-1,n+1}
\right)\\ & = -in f_{m+1,n-1} -im f_{m-1,n+1}.
\end{align*}
The proposition follows since $W_m$ obviously is closed under all 
three basis operators.
\end{proof}
\begin{remark}
Using the three basis operators given above it is easy to see that $W_m$ 
is irreducible.
\end{remark}

Instead of choosing the basis $\{f_{m-n,n}\}_{n=0}^m$ for $W_m$ 
it is sometimes more convenient to use the basis of eigenfunctions of 
$d \omega \left( \begin {array}{ccc} i & 0\\
\noalign{\medskip} 0 & -i\end {array} \right).$ Because of the symmetry in the 
commutator relations of the basis elements of $\mathfrak{su}(2,\C),$ these 
eigenfunctions have the same set of eigenvalues as $\{f_{m-n,n}\}_{n=0}^m.$ Call 
this new basis $\{b_{m,n}\},$ where $n=-m,-m+2,...,m$ and $b_{m,n}(re^{i\theta})
=e^{in\theta}b_{m,n}(r).$ The elements of the basis is determined by the 
relations above up to multiplication by a constant, choosing these 
constants correctly we get:

\begin{prop}
Let 
$$L_n^{(\alpha )}(x)= \frac{x^{-\alpha}e^x}{n!}\frac{d^n}{dx^n}\left(x^{n+\alpha}
e^{-x}\right)$$
be the Laguerre polynomials. (See \cite{Sze})
We have that
$$b_{m,n}(re^{i\theta})=e^{in\theta}r^{|n|}L_{(m-|n|)/2}^{(|n|)}(2\pi r^2)
e^{-\pi r^2}.$$
\end{prop}
\begin{proof}
We assume $n\ge 0,$ the argument is same as for $n<0.$ 
Since $b_{m,n}\in W_m,$ we see that $b_{m,n}$ is on the form 
$c(z, \bar z)e^{-\pi |z|^2},$ where $c$ is a polynomial of degree $m.$ That
$b_{m,n}(re^{i\theta})=e^{in\theta}b_{m,n}(r)$ 
means that $c(z,\bar z)$ only consists of terms on 
the form $z^a\bar z^b,$ where $a-b=n.$ In particular we must have that 
$b_{m,n}(re^{i\theta})=e^{in\theta}r^{n}q_{m,n}(2\pi r^2)e^{-\pi r^2},$
where $q_{m,n}$ is a polynomial of degree $(m-n)/2.$
Since the subspaces $W_m$ are orthogonal to each other, for $m\not =m'$ we have
\begin{align*}
0 & =\langle\overline{b_{m,n}},b_{m',n}\rangle=2\pi \int_0^\infty r^{n}
\overline{q_{m,n}(2\pi r^2)}
e^{-\pi r^2}r^{n}q_{m',n}(2\pi r^2)e^{-\pi r^2}r dr\\
&=\frac{1}{2(2\pi)^n}\int_0^\infty \overline{q_{m,n}(x)}q_{m',n}(x)x^{n}e^{-x}dx.
\end{align*}
This proves that $q_{m,n}(x)=L_{(m-n)/2}^{(n)}(x)$ if we normalize correctly. 
\end{proof}

\section{Properties of the local Tate $\zeta$-function}
\begin{define}
We define the local Tate $\zeta$-function
\begin{align*}
\zeta(s,\nu,f)=\int_{\C^{\times}} f(z) \nu(z) |z|^{2s-2}dz
\end{align*}
for all characters $\nu$ of $\C^{\times}$ and $f\in S(\C).$ 
\end{define}
\begin{remark}
This is the local $\zeta$-functions defined in the introduction specialized 
to the case where the local field is $\C$.
\end{remark}
All characters of 
$\C^{\times}$ can be written using polar coordinates 
in the form $\nu(r,\theta )=
r^{i\alpha}e^{ik\theta}$ with 
$k\in \mathbb{Z}.$ Since 
$\zeta(s,r^{i\alpha}e^{ik\theta},f)=\zeta(s+i\alpha /2,e^{ik\theta},f),$
the real part of the zeros of $\zeta$ does not depend on $\alpha.$ 
Hence our attention will be drawn to the following object:
\begin{define}
Let $k,m\in \mathbb{N},$ $\nu_k=e^{ik\theta}$ and $g_k=r^{2s-2}\nu_k.$ 
We set
\begin{align*} 
\zmk(s)=\langle f_{m,0},g_k\rangle=\zeta(s,\nu_k,f_{m,0}).
\end{align*}
\end{define}

In order for Theorem \ref{lrh} to be true it is essential that 
all elements in the invariant subspaces define the same $\zeta$-function 
$\zmk,$ up to multiplication by a constant. 
That this really is the case is shown in 
the next proposition. 

\begin{prop}
\label{unikhet}
If $f\in W_m$ then $\zeta(s,\nu_k,f)=c_{f,k}\cdot\zmk(s),$ where $c_{f,k}$ is a 
constant not depending on $s.$
\end{prop}

\begin{proof}
Let $f=\sum_{j=0}^mc_{2j-m}b_{m,2j-m}.$ For $(m-k)/2\in \mathbb{N}$ we see that
\begin{align*}
\zeta(s,\nu_k,f)&=\sum_{j=0}^mc_{2j-m}\zeta(s,\nu_k,b_{m,2j-m})\\
&=\sum_{j=0}^mc_{2j-m}\int_0^\infty\int_0^{2\pi}e^{i(2j-m)\theta}b_{m,2j-m}(r)
r^{2s-1}e^{ik\theta}d\theta dr\\
&=c_k\zeta(s,\nu_k,b_{m,k}),
\end{align*}
other $m$ give $\zmk(s)\equiv 0.$
\end{proof}

\begin{lemma}
\label{polynom}
If $(m-k)/2\in \mathbb{N}$ we have that 
$$\zmk(s)=\Gamma\left(s+\frac{k}{2}\right)\pi^{1-s}p_m^{(k)}(s),$$
where $p_m^{(k)}(s)$ is a real polynomial of degree $(m-k)/2.$ 
Otherwise $\zmk(s)\equiv 0.$
\end{lemma}
\begin{proof}
Since $H_m$ is odd if $m$ is odd and even if $m$ is even, the trigonometric 
identities \cite{grad}
\begin{align*}
\cos ^{2n}\theta=\frac{1}{2^{2n}}
\left( \begin {array}{ccc} 2n\\ \noalign{\medskip} n\end {array} \right)
+\frac{1}{2^{2n-1}}\sum_{j=1}^n
\left( \begin {array}{ccc} 2n\\ \noalign{\medskip} n-j\end {array} \right)
\cos (2j\theta)
\end{align*}
and
\begin{align*}
\cos ^{2n-1}\theta=\frac{1}{2^{2n-2}}
\sum_{j=1}^n
\left( \begin {array}{ccc} 2n-1\\ \noalign{\medskip} n-j\end {array} \right)
\cos ((2j-1)\theta),
\end{align*}
can be used to write
$$H_m(\sqrt{2\pi}r\cos \theta)=\sum_{j=0}^{\left[m/2\right]}
r^{m-2j} a_j(r^2)\cos((m-2j)\theta)$$
for some real polynomials $a_j(r)$ with $\deg a_j=j.$ 
This implies that if $(m-k)/2\not\in \mathbb{N}$ then 
$\zmk(s)\equiv 0$ and if $(m-k)/2\in \mathbb{N}$ we have
\begin{align*}
\zmk(s)&=\int_0^\infty\int_0^{2\pi}H_m(\sqrt{2\pi}
r\cos\theta)e^{-\pi r^2}r^{2s-1}
e^{ik\theta}d\theta dr\\
&=2\pi\int_0^\infty r^k a_{\frac{m-k}{2}}(r^2)r^{2s-1}e^{-\pi r^2}dr=
\pi\sum_{j=0}^{\frac{m-k}{2}} b_j \int_0^\infty r^{2s-1+k+2j}e^{-\pi r^2}dr\\
&=\pi\sum_{j=0}^{(m-k)/2} b_j\frac{1}{2\pi^{s+j+k/2}}
\Gamma\left(s+j+\frac{k}{2}\right)\\
&=\sum_{j=0}^{(m-k)/2} b_j\frac{1}{2\pi^{s+j+k/2-1}}
\left(s+j+\frac{k}{2}-1\right)...\left(s+\frac{k}{2}\right)
\Gamma\left(s+\frac{k}{2}\right)\\
&=\Gamma\left(s+\frac{k}{2}\right)\pi^{1-s}p_m^{(k)}(s),
\end{align*}
where $p_m^{(k)}(s)$ is a real polynomial of degree $(m-k)/2.$
\end{proof}

\begin{remark}
Theorem \ref{lrh} implies that 
$p_m^{(k)}(1-s)=(-1)^{\frac{m-k}{2}}p_m^{(k)}(s)$
so $\zmk$ fulfills a functional equation much like the functional equation for 
the Riemann $\zeta$-function.
\end{remark}

\begin{lemma}
\label{fekv}
$\zmk(s)$ admits the functional equation
\begin{align*}
(m+1)\zmk (s) = \pi \zmk(s+1)-\frac{1}{\pi}\left(s+\frac{k}{2}-1\right)
\left(s-\frac{k}{2}-1\right)\zmk(s-1).
\end{align*}
\end{lemma}

\begin{proof}
Since we have that
\begin{align*}
\Delta g_k(s)=\left(\frac{\partial^2}{\partial r^2}+\frac{1}{r}\frac{\partial}
{\partial r}+\frac{1}{r^2}\frac{\partial^2}
{\partial \theta^2}\right)r^{2s-2}e^{i\theta k}=
\left((2s-2)^2-k^2\right)g_k(s-1)
\end{align*}
and 
\begin{align*}
\left(-\frac{1}{4\pi}\Delta +\pi\left(x^2+y^2\right)\right)f_{m,0}=
(m+1)f_{m,0},
\end{align*}
we immediately get
\begin{align*}
(m+1)\zmk(s)&=\langle (m+1)f_{m,0},g_k(s)\rangle
=\left\langle\left(-\frac{1}{4\pi}\Delta+\pi 
\left(x^2+y^2\right)\right)f_{m,0},g_k(s)\right\rangle\\&=
\left\langle f_{m,0},\left(-\frac{1}{4\pi}\Delta+\pi 
\left(x^2+y^2\right)\right) g_k(s)\right\rangle\\&=
\left\langle f_{m,0},-\frac{1}{4\pi}
\left((2s-2)^2-k^2\right)g_k(s-1)+\pi g_k(s+1)\right\rangle\\&=
-\frac{1}{\pi}\left(s+\frac{k}{2}-1\right)
\left(s-\frac{k}{2}-1\right)\zmk(s-1)+\pi \zmk(s+1).
\end{align*}
\end{proof}

From \cite{lrh} we have the following lemma: 

\begin{lemma}
\label{metod2}
Let $q(s)$ be a polynomial, and assume that the zeros of $q(s)$ lie in
the closed strip $\{s;\Re(s)\in [-c,c]\}$ with $c>0.$ 
Then if $a,b>0,$ the zeros of 
\begin{align*}
r(s)=(s+a)q(s+b)-(s-a)q(s-b)
\end{align*}
lie in the open strip  $\{s;\Re(s)\in (-c,c)\}.$
\end{lemma}

\begin{remark}
The lemma is proved for $b=2$ but this does not change the proof.
\end{remark}

\begin{proof} [Proof of Theorem \ref{lrh}]
We only need to show that $p_m^{(k)}(s)$ has all its zeros on $\Re (s)=1/2.$
Letting $q_m^{(k)}(s)=p_m^{(k)}(s+1/2)$ and inserting this in 
Lemma~\ref{fekv} we get
\begin{align*}
(m+1)\Gamma \left(s+\frac{k+1}{2}\right)\pi^{\frac{1}{2}-s}q_m^{(k)}(s)=
\pi\Gamma \left(s+1+\frac{k+1}{2}\right)\pi^{-\frac{1}{2}-s}q_m^{(k)}(s+1)\\
-\frac{1}{\pi}\left(s+\frac{k-1}{2}\right)\left(s-\frac{k+1}{2}\right)
\Gamma \left(s-1+\frac{k+1}{2}\right)\pi^{\frac{3}{2}-s}q_m^{(k)}(s-1).
\end{align*}
Simplifying this gives
\begin{align*}
(m+1)q_m^{(k)}(s)=\left(s+\frac{k+1}{2}\right)q_m^{(k)}(s+1)
-\left(s-\frac{k+1}{2}\right)q_m^{(k)}(s-1).
\end{align*}
The claim now follows from Lemma~\ref{metod2}.
\end{proof}

We could also prove Theorem \ref{lrh} in a different way by using the 
following well-known theorem:

\begin{thm}
\label{metod1}
Let $\{p_n\}_{n=0}^\infty$ be a sequence of polynomials such that the degree 
of $p_n$ is $n$ and the polynomials are orthogonal with respect to some Borel 
measure $\mu$ on $\reals .$ Then $p_n$ have $n$ distinct real roots.
\end{thm}
\begin{proof}
The theorem is obviously true for $n=0.$ Assume that $p_k$ has $k$ distinct 
roots for $k<n.$ Without loss of generality we assume that all polynomials have 
one as their leading coefficient. Then $p_k$ is real for $k<n$ and 
$p_n=f_n+ig_n,$ where $g_n$ has degree less than $n.$ Moreover, 
\begin{align*}
0=(p_n,p_k)=(f_n,p_k)-i(g_n,p_k)
\end{align*}
for $k<n,$ hence $(g_n,p_k)=0.$ But the degree of $g_n$ is less than $n$ so we 
must have $g_n\equiv 0.$ Thus $p_n$ is real. If $p_n$ does not have $n$ distinct 
real roots then it could be written as 
$p_n(x)=(x-\alpha)(x-\bar\alpha)q(x)=|x-\alpha|^2q(x)$ for $x\in\reals .$ 
Since the degree of $q$ is less than $n$ we must have $(p_n,q)=0,$ but on the 
other hand we have that $p_n(x)q(x)\ge 0$ for all $x.$ This is a contradiction, 
hence $p_n$ has $n$ distinct real roots.
\end{proof}

\begin{prop}
\label{orto}
The polynomials $p_m^{(k)}(1/2+it)$ are orthogonal with respect to the measure 
$|\Gamma((k+1)/2+it)|^2dt,$ where $dt$ is the Lebesgue measure on $\reals .$
\end{prop}

\begin{proof}
As we have noticed before the functions $b_{m,n}$ are orthogonal, hence for 
$m\not =m'$ we have
\begin{align*}
0 & =\langle\overline{b_{m,n}},b_{m',n}\rangle=
2\pi\int_0^\infty \overline{b_{m,n}(r)}b_{m',n}(r)rdr\\ 
& =2\pi\int_{-\infty}^\infty\overline{b_{m,n}(e^u)e^u}b_{m',n}(e^u)e^udu.
\end{align*}
Using Plancherel's formula it follows that $2\pi \mathfrak{F} 
(b_{m,n}(e^u)e^u)(-2t)$ is an orthogonal sequence ($\mathfrak{F}$ denotes 
the ordinary Fourier transform) and this is just
\begin{align*}
2\pi\mathfrak{F}(b_{m,k}(e^u)e^u)(-2t)& = 2\pi\int_0^\infty 
b_{m,-k}(r)r^{i2t}dr\\
&=\int_0^\infty \int_0^{2\pi}b_{m,-k}(re^{i\theta})e^{ik\theta}r^{i2t}d\theta dr
=c_{m,k}\zmk (1/2+it)\\
&=c_{m,k}\Gamma\left(\frac{k+1}{2}+it\right)\pi^{1/2-it}p_m^{(k)}(1/2+it).
\end{align*}
\end{proof}
\begin{remark}
Theorem \ref{lrh} follows immediately if we combine Theorem~\ref{metod1} with 
Proposition~\ref{orto}.
\end{remark}


\bibliography{ref}
\bibliographystyle{plain}

\end{document}